\documentclass{amsart}
\address{Alan T. Huckleberry, Ruhr-Universit\"at Bochum, Fakult\"at f\"ur Mathematik,
44780 Bochum, Germany, {\tt ahuck@cplx.ruhr-uni-bochum.de}}
\author{Alan T. Huckleberry, Stefan Kebekus and Thomas Peternell}
\address{Stefan Kebekus, Lehrstuhl Mathematik VIII, Universit\"at Bayreuth,
95440 Bayreuth, Germany, {\tt stefan.kebekus@uni-bayreuth.de}}
\address{Thomas Peternell, Lehrstuhl Mathematik I, Universit\"at Bayreuth,
95440 Bayreuth, Germany, {\tt peternel@btm8x1.mat.uni-bayreuth.de}}

\usepackage{amssymb}
\usepackage{amscd}
\usepackage[arrow,matrix,curve,ps]{xy}
\xyoption{dvips}

\theoremstyle{plain}    
\newtheorem{thm}{Theorem}[section]
\numberwithin{equation}{section} 
\numberwithin{figure}{section} 
\theoremstyle{plain}    
\newtheorem{cor}[thm]{Corollary} 
\theoremstyle{plain}    
\newtheorem{lem}[thm]{Lemma} 
\theoremstyle{plain}    
\newtheorem{prop}[thm]{Proposition} 
\theoremstyle{definition}

\theoremstyle{remark}    
\newtheorem{notation}[thm]{Notation}

\begin{document}

\title{Group Actions on $S^6$ and complex structures on $\mathbb P_3$}
\begin{abstract}
It is proved that if $S^6$ possesses an integrable complex structure,
then there exists a $1$-dimensional family of pairwise different
exotic complex structures on $\mathbb P_3(\mathbb C)$.  This follows
immediately from the main result of the paper: $S^6$ is not the
underlying differentiable manifold of an almost homogeneous complex
manifold $X$.  Via elementary Lie theoretic techniques this is reduced
to ruling out the possibility of a $\mathbb C^*$-action on a certain
non-normal surface $E\subset X$.  A contradiction is reached by
analyzing combinatorial aspects of the non-normal locus $N$ of $E$ and
its preimage $\hat N$ in the normalization $\hat E$.
\end{abstract}

\date{\today}

\maketitle
\tableofcontents

\section{Introduction}

This note is motivated by the following classical problem: Is there a
complex structure on the 6-sphere $S^{6}$, i.e. is $S^{6}$ a complex
manifold?

It has been known known since decades that all other spheres $S^{2n}$
do not even admit an almost complex structure. On the other hand $S^6$
admits many almost complex structures. It is generally believed that
none of them is integrable.

Suppose $S^{6}$ has a complex structure $X$. Then by \cite{CDP} every
meromorphic function on $X$ is constant. Moreover $X$ is not K\"ahler,
since $b_2(X)=0$.  Therefore the problem is quite inaccessible by
standard methods of complex geometry.

In this paper we prove

\begin{thm}\label{mainThm}
$X$ is not almost homogeneous. In other words, the automorphism group
$Aut_{\mathcal O}(X)$ does not have an open orbit.
\end{thm}

This is related as follows to the question of existence of complex
structures on the underlying differentiable manifold of $\mathbb
P_3(\mathbb C)$: As above assume $X=S^6$ has the structure of a
complex manifold.  For $p\in X$ let $\pi _p:X_p\to X$ denote the blow
up of $X$ at $p$ with $\pi _p^{-1}(p)=:E_p$.  Of course $E_p=\mathbb
P_2(\mathbb C)$.  Since sufficiently small neighborhoods of a
hyperplane in $\mathbb P_n(\mathbb C)$ are differentiably identifiable
with neighborhoods of a blown up point, it follows that $X_p$ is
diffeomorphic to $\mathbb P_3(\mathbb C)$.

Suppose that for given points $p,q\in X$ there exists a biholomorphic
mapping $\psi :X_p\to X_q$.  Note that, since $\psi (E_p)$ generates
the cohomology $H^2(X_q,\mathbb Z)$, it follows that $\psi (E_p)\cap
E_q\ne \emptyset$.  If $C:=\psi (E_p)\cap E_q\ne E_q$, then $\pi
_q\vert \psi(E_p)$ would be a modification which maps the curve
$C\subset \psi(E_p)\cong \mathbb P_2(\mathbb C)$ to a point.  Since
this is impossible, $\psi(E_p)=E_q$ and $\psi $ induces an
automorphism $g_\psi :X\to X$ with $g_\psi (p)=q$.

Let ${\mathcal F}$ denote the orbit space $X/Aut_{\mathcal O}(X)$.
For $\xi,\eta \in {\mathcal F}$ with $\xi \ne \eta $ and $p\in \xi ,
q\in \eta $ representatives, it follows that $X_p$ and $X_q$ are not
biholomorphically equivalent.

Of course ${\mathcal F}$ may very well be non-Hausdorff and, from
certain points of view, an unreasonable parameter space, but by abuse
of language we nevertheless refer to it as a family of complex
structures on $\mathbb P_3(\mathbb C)$.  Since $Aut_{\mathcal O}(X)$
does not have an open orbit, its generic orbit in $X$ is at least
$1$-codimensional.  In particular, for general $p$, the semi-universal
deformation space of $X_p$ is at least 1-dimensional. In this sense we
have the following consequence of theorem~\ref{mainThm}.

\begin{cor}
If $S^{6}$ admits a complex structure, then there is a 1-dimensional
family of complex structures on $\mathbb P_3$.
\end{cor}

\begin{cor}\label{Intro:Cor2}
Let $X$ be a complex structure on $S^{6}$. Then $X$ carries at most
two linearly independent holomorphic vector fields: $h^0(TX)\leq 2$.
\end{cor}

\begin{proof}
If the generic $Aut_{\mathcal O}(X)$-orbit is $2$-dimensional, then
there are globally defined holomorphic vector fields $V_1$ and $V_2$
such that $U:=\{ p\in X: V_1\wedge V_2(p)\ne 0\} $ is a dense, Zariski
open subset.  If $V_3$ is any holomorphic field on $X$, then, since
$V_1\wedge V_2 \wedge V_3\equiv 0$, there exist $m_1,m_2\in {\mathcal
O}(U)$ with $V_3=m_1V_1+m_2V_2$ on $U$.  By explicit computation in
local coordinates, on verifies that these coefficients are in fact
meromorphic on $X$ and the relation between the fields extends to
$X$. On the other hand, it has been proven in \cite{CDP} that $X$
possesses only the constant meromorphic functions. Thus, it follows in
this case that $dim_{\mathbb C}\Gamma (X,TX)=2$. The other cases,
i.e., those where the generic orbit dimension is smaller, are handled
analogously.
\end{proof}
By semicontinuity we obtain from Corollary \ref{Intro:Cor2}:

\begin{cor}
The inequality $h^0(TX\otimes L)\leq 2$ holds for generic $L\in
Pic^0(X)$.
\end{cor}
Our theorem should be seen in a more general context or rather program
for investigating complex structures on $S^{6}$. Namely we would like
to prove that the tangent bundle $TX$ is stable with respect to a
Gauduchon metric. It would then follow that $X$ carries a
Hermite-Einstein connection (Li-Yau). At this point one could employ
powerful analytic tools to investigate the problem further. In order
to prove the stability it would seem necessary to verify the
statements:
\begin{description}
\item [(A)]$H^0(X,TX\otimes L)=0$ for all $L\in Pic^0(X)$;
\item [(B)]$H^0(X,\Omega _{X}^1\otimes L)=0$ for all $L\in Pic^0(X)$.
\end{description}
Hence our theorem is the first approach to (A). The next step should
be to investigate group actions in general. We feel that the study of
group actions on highly non-algebraic manifolds is of independent
interest and hope that the methods which we develop in this paper are
useful in a broader context.

\section{Setup and general results}\label{sect:setup}

In this section we gather results which will be used throughout the
paper and give the general setup.

\subsection{Setup and outline of proof}

Let $X$ denote a complex structure on $S^{6}$. The entire paper is
devoted to proving that $X$ is not almost homogeneous.  This is proved
by assuming the contrary, i.e., there exists $x_0\in X$ such that
$G.x_0=:\Omega $ is open, and deriving a contradiction.

The structure of the proof can be outlined as follows: Using the fact
that $X$ has only finitely many analytic hypersurfaces in connection
with with topology and Lie theory of the situation at hand, it is
shown in sections \ref{sect:semisimple} that if $G$ exists, then it
must be solvable. In section \ref{sect:methods}, we give a number of
methods which are applied in section \ref{sect:elim} in order to rule
out this case, too.

Two of the methods involve a lengthy proof which we have preferred to
give separately in section \ref{CCAA:Section} and sections
\ref{StartC*}--\ref{sect:fixed_curve}.

The technical heart of this work lies in sections
\ref{StartC*}--\ref{sect:fixed_curve} where we rule out the following
situation: we suppose that there exists of a subgroup $\mathbb
C^*<Aut_{\mathcal O}(X)$ and that $E:=X\setminus \Omega $ is an
irreducible, non-normal, rational surface where $\mathbb C^*$ acts as
an algebraic transformation group.

The non-normal locus $N\subset E$ and its preimage $\tilde N$ in the
normalization $\tilde E$ play a key role in the remainder of the
proof.  It follows from the Betti number information on $X$ that $N$
and $\tilde N$ are connected and have the same first Betti numbers. An
analysis of the $\mathbb C^*$-action on $E$ shows however that in fact
$b_1(N)=b_1(\tilde N)+1$ (see sections
\ref{sect:discrete_fixed_points}--\ref{sect:fixed_curve}).

\subsection{Meromorphic functions, discrete isotropy, and the
dimension of $E$}

Since there are no non-constant meromorphic function on $X$ by
\cite{CDP}, we have that $\dim G=h^0(TX)=3$.

\begin{prop}
The $G$-isotropy $\Gamma$ at a point $x_0\in \Omega$ is
discrete. Furthermore, $E:=X\setminus \Omega$ is non-empty and
1-codimensional in $X$.
\end{prop}
\begin{proof}
Since $dim\,G=3$, it is clear that $\Gamma $ is discrete. In
particular, if $G$ would act transitively, then, contrary to
assumption, $\pi _1(X)\cong \Gamma $ would not be trivial. Thus $E\ne
\emptyset $.  If the vector fields $X_1,X_2,X_3$ form a basis of
$\Gamma (X,TX)$, then
$$
E=\{ X_1\wedge X_2\wedge X_3 =0\} .
$$ 
In particular, $-K_X=\mathcal{O}_X(\sum\lambda_i E_i)$, where $E_{i}$
are the irreducible components of $E$ and $\lambda _{i}>0$.
\end{proof}

Let $G=R.S$ be the Levi-Malcev decomposition, so that $R$ is
the radical of $G$, i.e. the maximal connected solvable normal
subgroup of $G$ and $S$ is semisimple; moreover $R\cap S$ is discrete.

Since a semi-simple complex Lie group has dimension at least 3,
we have only two cases, namely that $G$ is semi-simple or solvable.

\subsection{Topology of $\Omega$ and $E$}

Since the topology of $X$ is well-known, the Betti-numbers of
$\Omega $ and $E$ are closely related.

\begin{notation}
If $Y$ is a complex space, set $b_{i}(Y):=h^{i}(Y;\mathbb Q)$.
\end{notation}

\begin{prop}\label{Top:EvsOmega}
For $q\in \{1,\ldots ,4\}$, the following equation of
Betti-numbers holds: $b_{q}(E)=b_{5-q}(\Omega )$.
\end{prop}

\begin{proof}
Recall from algebraic topology that there is an exact cohomology
sequence associated to the pair $(X ,E)$:
$$
\ldots \to H^{q}(X;\mathbb Q)\to H^{q}(E;\mathbb Q)\to H^{q+1}
(X,E;\mathbb Q)\to H^{q+1}(X,\mathbb Q)\to \ldots 
$$
Since $X$ is homeomorphic to the 6-sphere, $b_{q}(X)=b_{q+1}(X)=0$ for
all numbers $1\leq q\leq 4$, so that $H^{q}(E;\mathbb Q)\cong
H^{q+1}(X,E;\mathbb Q)$.  An application of the Alexander
duality theorem yields $H^{q+1}(X,E;\mathbb Q)\cong H_{5-q}(\Omega
;\mathbb Q)$, hence the claim.
\end{proof}

\subsection{Fixed points of reductive groups}

Many of our arguments involve linearization of group actions at 
fixed points. We recall the theorem on faithful linearization:

\begin{thm}\label{linearization}
Assume that a reductive complex Lie group $H$ acts
holomorphically on a complex manifold $M$, and assume that $x\in M$ is
an $H$-fixed point. For $h\in H$, let $T(h):T_{x}\to T_{x}$ be the
tangential map. Then there exist neighborhoods $U$ of $x$, and $V$ of
$0\in T_{x}M$ and an isomorphism $\phi :U\to V$ such that $\phi \circ
h=T(h)\circ \phi $ for all $h$ in a given maximal compact subgroup of
$H$.

Furthermore, if $W$ is a neighborhood of the maximal compact subgroup
and $U'\subset U$ open so that $WU'\subset U$, then $(T(w)\circ \phi
)(x)=(\phi \circ w)(x)$ for all $x\in U'$.
\end{thm}
In this setting we call $U$ a \emph{linearizing neighborhood} of $x$.
See \cite{Huc90} or \cite[p. 11f]{H-Oe} information about
linearization.

The fixed point set of a reductive group also possesses certain
topological properties. In our special case this implies:

\begin{prop}\label{Fix.of.C*}
Suppose that $Aut(X)$ contains a subgroup $I\cong \mathbb C^{*}$.  Let
$F$ be the fixed point set of $I$. Then either $F\cong \mathbb P_1$,
or $F$ consists of two disjoint points.
\end{prop}

\begin{proof}
As a first point, note that it follows from the linearization theorem
that $F$ is smooth; in particular, its irreducible components are
disjoint.  Furthermore, $\chi _{Top}(F)=\chi _{Top}(X)$ (see \cite{KP85} for
a proof in the algebraic setting which carries over immediately to,
e.g., compact complex spaces.)

In our situation, for $p\in F$ the group $I$ stabilizes the complement
$X\setminus \{ p\} \cong S^6$ and thus $F$ is acyclic (see
\cite{Br72}).  Consequently, $F$ consists of two points or it is
irreducible. Thus it remains to show that $F$ is at most
$1$-dimensional.

Suppose $F$ is $2$-dimensional. Since $H^4(X,\mathbb Z)=0$, it follows
that $c_2(TX|_F)=0$ and, since $H^2(X,\mathbb Z)=0$, the normal
bundle $N_{F,X}$ is topologically trivial.  Consequently
$$
0=c_2(TX|_F)=c_2(F)+N_{F,X}.K_F=\chi _{Top}(F)
$$
which is contrary to $\chi _{Top}(F)=2$.
\end{proof}

\section{The case where $G$ is semisimple}\label{sect:semisimple}

In this section we treat the case that $G$ is semisimple. Since $\dim
G=3$, it follows that $G\cong SL_2$. The most basic property of
almost transitive $SL_2$-actions on threefolds is:

\begin{lem}
\label{SL_2:no.fixed.points.}The $G$-action on $X$ does not have
a fixed point. In particular, if $\tilde{E}_{i}$ is the normalization of
$E_{i}$, then $\tilde{E}_{i}$ is smooth.
\end{lem}
\begin{proof}
Assume that $x\in X$ was a $G$-fixed point. Linearize the $G$-action
at $x$ and recall from the representation theory of $SL_2$ that no
3-dimensional $SL_2$-representation space is
$SL_2$-almost homogeneous (see \cite[lemma 1.12]{Mu83} for a more
detailed proof).  A contradiction.
\end{proof}

\begin{lem}
If $E_{i}\subset E$ is an irreducible component, then $G$ acts almost
transitively on $E_{i}$. In particular, the normalization $\tilde{E}_{i}$
is either rational or a Hopf-surface.
\end{lem}

\begin{proof}
If $G$ did not have an open orbit in $E_i$, then by
lemma~\ref{SL_2:no.fixed.points.} all $G$-orbits would be
$1$-dimensional. Thus the normalization $\tilde E$ would be the
product $\mathbb P_1(\mathbb C)\times C$, where $G$ operates
transitively on the first factor and $C$ is a smooth curve.  In
particular it follows that a maximal torus $T\cong \mathbb C^*$ would
have two disjoint copies of $C$ as a fixed point set in $E_i$.  Note
that, since the normalization map $\tilde E_i \to E_i$ is equivariant
with respect to $G\cong SL_2$, the $T$-fixed point set in $E_i$ is the
disjoint union of two curves, contrary to proposition~\ref{Fix.of.C*}.

The ``In particular$\ldots $'' clause of the statement results from
the classification of the almost homogeneous surfaces ---see
\cite[p. 92]{H-Oe} or \cite{Pot69}.
\end{proof}

Now we exclude both possibilities:

\begin{prop}\label{G_not_semisimple}
The group $G$ is not semisimple.
\end{prop}
\begin{proof}
It follows from the above lemma that the normalization $\tilde E_i$ of
a component of $E$ is an almost homogeneous Hirzebruch surface,
$\mathbb P_2(\mathbb C)$ equipped with either the defining
representation of $SO_3(\mathbb C)$ or the representation of $SL_2$
with a fixed point or a homogeneous Hopf surface (see \cite{H-Oe} or
\cite{Pot69}). Consequently, a maximal torus $T\cong C^*<G$ has only
isolated fixed points in $X$ and, if $E_i$ were rational it would
already have $3$ or $4$ fixed points in $E_i$ alone.

Thus we may assume that every such component is a homogeneous Hopf
surface. But this is also not possible, because such a surface is a
homogeneous space $G/H$, where $H^0$ is unipotent, i.e., $T$ has no
fixed points.
\end{proof}

\section{Main methods for the elimination of solvable groups}
\label{sect:methods}

We begin by presenting several methods which involve the normalizer of
subgroups of isotropy groups. Recall that $G$ is a connected,
simply-connected complex Lie group acting almost transitively
on $X$ with an open orbit $\Omega = G.x_0$. By
proposition~\ref{G_not_semisimple} we may assume that $G$ is solvable
and by the remarks in section \ref{sect:setup} that the isotropy
$\Gamma:=G_{x_0}$ is discrete.

\subsection{The normalizer arguments}

Throughout the paper if $H$ is a subgroup of $G$, then $N(H)$ denotes
its normalizer in $G$ and $N(H )^0$ its identity component. If $H$ is
discrete, it follows that $N(H)^0=Z(H)^0$, where $Z(H)$ denotes its
centralizer.

\subsubsection{The 2-dimensional normalizer argument}
Note that if $H<\Gamma $ is
not normal in $G$, then it acts non-trivially on $\Omega $.

\begin{prop}[2-dimensional normalizer argument]\label{2dNormArg}
If $ H<\Gamma $ is an arbitrary subgroup, then $\dim N(H)\not =2$.
\end{prop}

\begin{proof}
Suppose not and note that the 2-dimensional orbit $F:=Z(H)^0.x_0$ is a
component of the set of $H$-fixed points in $\Omega$. Since the full
set $X^H$ of $H$-fixed points is closed, $F$ is Zariski open in
it's closure $\overline F$ which is 1-dimensional in $X$.

Observe that $\{g\overline F|g\in G\}$ is an infinite set of
hypersurfaces and consequently there exists a non-constant meromorphic
function on $X$ (see \cite[thm. 1]{Kras75} \footnote {See \cite{FF79}
for a more general result.}). However, such a function does not exist
(see \cite{CDP}).
\end{proof}

\subsubsection{The 1-dimensional normalizer argument}

\begin{prop}[1-dimensional normalizer argument]\label{1dNormArg}
If $H<\Gamma$ is any subgroup and $\dim N(H)=1$, then $N(H)^0\cap
\Gamma$ is a lattice of rank 2.
\end{prop}
\begin{proof}
The 1-dimensional group $Z:=Z(H)^0$ acts transitively on the fixed
point components $\Omega^H\subset \Omega$ which are of course Zariski
open in their closures. Set $C:=Z.x_0$ and note that if $C\not =
\overline C$, then $\overline C.E>0$, contrary to $h^2(X;\mathbb
Z)=0$. Thus the orbit $C$ is a compact 1-dimensional torus.
\end{proof}

\subsection{Arguments involving reductive subgroups of $Aut(X)$}

In many cases we are able to rule out the existence of certain
subgroups of $Aut(X)$ under additional assumptions on the topology of
$\Omega$. For convenience, use the following

\begin{notation}
If $Y$ and $Z$ are topological spaces, say that $Y$ has
``Betti-type $Z$'' if all the Betti-numbers of $Y$ and
$Z$ agree.
\end{notation}

\subsubsection{The $\mathbb C^{*}\times \mathbb C^{*}$-action argument}

\begin{prop}\label{CCAA:prop}
If $E$ is connected and has at least two irreducible components, then
there does not exist an (effective) action of $(\mathbb C^*)^2$ on $X$.
\end{prop}
The proof will be given in section \ref{CCAA:Section}.

\subsubsection{The torus-action argument}

\begin{prop}[torus-action argument]\label{TorusActArg}
The automorphism group of $X$ does not contain a compact torus. In
particular, if $\Gamma '<\Gamma $ is an arbitrary subgroup which is
normal in $G$ and contained in a 1-dimensional Abelian subgroup $A$,
then $rank(\Gamma ')=1$.
\end{prop}
\begin{proof}
Assume that $T<Aut(X)$ is a 1-dimensional complex torus. Since $\chi
_{Top}(X)=2$, the Lefschetz fixed point formula shows that
every vector field must have zeros. In particular, $T$ must have a
fixed point in $X$. On the other hand, all representations of $T$ are
trivial, a contradiction the theorem on faithful linearization!
\end{proof}

\subsubsection{The $\mathbb C^{*}$-action argument}

\begin{prop}\label{S1Proposition}
If the open orbit $\Omega \subset X$ has the same Betti numbers
as the circle $S^1$, then the automorphism group of $X$ does not
contain a subgroup which is isomorphic to $\mathbb C^{*}$.
\end{prop}
The proof of the preceding proposition turns out to be a rather
involved matter.  We give it in sections~\ref{StartC*}--\ref{EndC*}.

The following corollary is useful when it comes to the exclusion of
certain subgroups $\Gamma $.

\begin{cor}[$\mathbb C^*$-action argument]\label{C*ActArg}
If the open orbit $\Omega \subset X$ has Betti-type $S^1$, and
if $\Gamma '<\Gamma $ is an arbitrary subgroup which is normal in $G$,
then $\Gamma '$ is not contained in any positive-dimensional Abelian
subgroup of $G$.
\end{cor}
\begin{proof}
Assume to the contrary: let $A<G$ be an Abelian, positive-dimensional
connected subgroup and consider the $A$-action on $X$, given by a
morphism $\rho :A\to Aut(X)$ of complex Lie groups.

Since it is normal, $\Gamma '$ acts trivially on $\Omega $, and thus
trivially on $X$. Consequence: the morphism $\rho $ factors through
the quotient $A/\Gamma '$, i.e. there exists an action of $A/\Gamma '$
on $X$.

Let $A'<A$ be a 1-dimensional subgroup such that $A'\cap \Gamma '\not
=\{0\}$.  Then $A'/(A'\cap \Gamma ')$ acts on $X$. By proposition
\ref{S1Proposition}, $A'/(A'\cap \Gamma ')$ must be a compact torus,
but we have seen in proposition \ref{TorusActArg} that no compact
torus acts non-trivially on $X$.
\end{proof}

\subsection{Topological observations}

Recall that $G$ is solvable and simply-connected, e.g.~it is a cell,
and $\Omega=G/\Gamma$, where $\Gamma$ is discrete. Thus the topology
of $\Omega$ is completely determined by $\Gamma$. It follows directly
from proposition \ref{Top:EvsOmega} that $\Gamma\not = \{e\}$.
 
\subsubsection{Restriction on the Betti-type of $\Omega$}

\begin{prop}[Betti-type argument]\label{TProposition}
The open orbit $\Omega \subset X$ does not have the Betti-type
of a $S^1\times S^1$.
\end{prop}
\begin{proof}
Assume to the contrary. Then, by proposition \ref{Top:EvsOmega}, we
have $b_{4}(E)=2$ and $b_3(E)=1$. Similar to the proof of proposition
\ref{Top:EvsOmega}, the sequence of the pair $(E,\Omega )$ and
Alexander duality give
$$
0\to \underbrace{H^0(X,E;\mathbb Q)}_{=H_{6}(\Omega ;\mathbb Q)=0}\to 
\underbrace{H^0(X;\mathbb Q)}_{=\mathbb Q}\to H^0(E;\mathbb Q)\to 
\underbrace{H^1(X,E;\mathbb Q)}_{=H_{5}(\Omega ;\mathbb Q)=0}\to \ldots 
$$
so that $b_0(E)=1$. In particular, $E$ is connected and has exactly
two irreducible components $E_1$ and $E_2$. Now $E_1|_{E_2}$ is a
Cartier divisor which is effective, non-zero, but homologeously
equivalent to zero, so that no desingularization of $E_2$ is
K\"ahler and vice versa.

By the Mayer-Vietoris sequence
$$
\ldots \to H^3(E)\to H^3(E_1)\oplus H^3(E_2)\to 
\underbrace{H^3(E_1\cap E_2)}_{=0}\to \ldots 
$$
we obtain $b_3(E_1)+b_3(E_2)\leq b_3(E)=1$. Hence we may assume that
$b_3(E_1)=0$.

Let $\nu :\tilde{E}_1\to E_1$be the normalization. Let $N\subset E$ be
the (set-theoretical) non-normal locus and $\tilde{N}:=\nu ^{-1}(N)$.

The sequence 
$$
\ldots \to H^{q}(E_1)\to H^{q}(\tilde{E}_1)\oplus H^{q}(N)\to 
H^{q}(\tilde{N})\to \ldots \quad q\geq 1
$$
yields $b_3(E_1)\geq b_3(\tilde{E}_1)$ so that
$b_3(\tilde{E}_1)=0$ (see \cite[prop. 3.A.7 on p. 98]{BK82} for an
explanation of the sequence).

Now consider the desingularization $\pi :\hat{E}_1\to
\tilde{E}_1$.  Since $\pi _{*}(\mathbb Z)$ is a skyscraper sheaf
with support only finitely many points,
$E_2^{p,q}:=H^{p}(\tilde{E}_1,R^{q}\pi _{*}(\mathbb Z))=0$ for all
$p>0$. Thus, $E_2^{p,q}\cong E_{\infty }^{p,q}$ in the Leray
spectral sequence, and since $R^3\pi _{*}(\mathbb Z)=0$ and
$H^3(\tilde{E}_1,\mathbb Z)=0$, we obtain $b_3(\hat{E}_1)=0$. Thus
$b_1(\hat{E})=0$, and the classification of surfaces yields that
$\hat{E}_1$ is K\"ahler\footnote{Actually it is easy to see
that $\kappa (\hat{E}_1)=-\infty $ so that we do not have to use the
fact that $K3$-surfaces are K\"ahler.}. This is a
contradiction!
\end{proof}

\section{Elimination of the solvable groups}
\label{sect:elim}

Our goal here is to eliminate the possibility that $X$ is almost
homogeneous with respect to the action of a solvable group. The proof
requires a bit of the special knowledge given by the classification
the of the simply connected 3-dimensional solvable Lie groups.

\subsection{Classification of the relevant Lie groups}

We now recall the classification mentioned above (see
e.g.~\cite{Jac62}). In this case $G$ is biholomorphic to $\mathbb C^3$
as a complex manifold and the group structure is isomorphic to one of
the following:

\begin{description}
\item [$\boldsymbol{G_0}$] this is the well-known Abelian group $\mathbb C^3$.
\item [$\boldsymbol{G_1}$] we could also denote this group by
$G_2(0)$. The multiplication is given as
$$
(a_1,b_1,c_1)(a_2,b_2,c_2)=(a_1+a_2,e^{a_1}b_2+b_1,c_1+c_2)
$$
\item [$\boldsymbol{G_2(\tau )}$] here $\tau $ is any complex number 
other than zero. The multiplication is
$$
(a_1,b_1,c_1)(a_2,b_2,c_2)=(a_1+a_2,e^{a_1}b_2+b_1,e^{\tau
a_1}c_2+c_1)
$$
\item [$\boldsymbol{G_3}$] Multiplication:
$$
(a_1,b_1,c_1)(a_2,b_2,c_2) =
(a_1+a_2,e^{a_1}b_2+b_1+a_1e^{a_1}c_2,e^{a_1}c_2+c_1)
$$
\item[$\boldsymbol{H_3}$] this is the Heisenberg group, where
the multiplication is
$$
(a_1,b_1,c_1)(a_2,b_2,c_2) =
(a_1+a_2,b_1+b_2+a_1c_2,c_1+c_2)
$$
\end{description}
For a detailed study of the discrete subgroups of such groups see
\cite[Sect. 1]{ES86}.

\subsection{The elimination}

Here we eliminate the possibility of an action of a solvable group by
utilizing a general strategy along with some knowledge of the
Lie groups which occur.

\begin{prop}\label{No_C3}
If $\Gamma$ is normal in $G$ and $A<G$ is any closed connected Abelian
subgroup, then $\Gamma\not\subset A$. In particular, the group $G$ is
not isomorphic to $G_0 \cong\mathbb C^3$.
\end{prop}
\begin{proof}
Assume to the contrary, i.e., that $\Gamma \subset A$.
Since $G/A$ is acyclic, $G/\Gamma $ has the same homotopy type as 
$A/\Gamma $, i.e., the same homotopy type as a real torus.

Our argument depends on $rank(\Gamma )$.
\begin{description}
\item [$\boldsymbol {rank(\Gamma )=1}$]Then $\Omega$ has the
Betti-type of the circle $S^1$.  Take a 1-dimensional subgroup $A'<G$,
$A'\cong \mathbb C$ with $\Gamma \subset A'$ and observe that
$A/\Gamma \cong \mathbb C^{*}$ acts non-trivially on $X$, contrary to
the $\mathbb C^*$-action argument~\ref{S1Proposition}.

\item[$\boldsymbol {rank(\Gamma )=2,3,4}$] Betti number considerations 
show that the assumptions of proposition~\ref{TorusActArg} are
satisfied. However, the assumption on $rank(\Gamma )$ together with
the torus-action argument allows us to construct a $(\mathbb
C^*)^2$-action coming from the group $A$.

\item[$\boldsymbol {rank(\Gamma )=5}$] Here $G=A$ is Abelian and 
$\Omega$ has two ends, i.e., proposition~\ref{TorusActArg} may not be
applied. However, the ineffectivity of the $G$-action on every
component of $E$ must contain a reductive group, i.e.~a torus or
$\mathbb C^*$, contrary to either the torus-action argument or the
characterization of the fixed point set of a $\mathbb C^*$ action.
\end{description}
This finishes the proof.
\end{proof}

A weaker statement holds if $\Gamma$ is not necessarily normal.

\begin{lem}\label{Abelian_lemma}
For all groups $A<G$ with $A\cong \mathbb C$ we have
$\Gamma\not\subset A$.
\end{lem}
\begin{proof}
Assume to the contrary. Again the Betti-types of $G/\Gamma$ and
of $A/\Gamma$ agree. We may assume that $\Gamma$ is not normal.

Since $A$ is contained in the normalizer of $\Gamma$, the
2-dimensional normalizer argument implies that $A=N_G(\Gamma)$. In
this setting the 1-dimensional normalizer argument yields that
$\Gamma\cong\mathbb Z^2$, contrary to the Betti-type argument.
\end{proof}

We will need the following technical lemma on centralizers of elements
in $G$.

\begin{lem}\label{centralizer_lemma}
If $G\cong G_1$, $G_2(\tau)$, $G_3$ or $H_3$ and $g\in G$, then $\dim
Z_G(g)\geq 1$ and one of the following holds:
\begin{enumerate}
\item $\dim Z_G(g)\geq 2$, or
\item $Z_G(g)^0 \subset (0,\mathbb C,\mathbb C)$, or
\item $Z_G(g)^0=Z_G(g)$, i.e.~$Z_G(g)$ is connected.
\end{enumerate}
Here $Z_G(g)$ denotes the set of elements commuting with $g$.
\end{lem}
\begin{proof}
The statement that $\dim Z_G(g)\geq 1$ can be checked directly. To
show that one of (1), (2) or (3) holds, for any number $a\in \mathbb
C$ define the set
$$
Z_a :=\{(b,c)\in\mathbb C^2|(a,b,c)\in G \text{ commutes with }g\}.
$$
Note that it is sufficient to show that for all $a$ either $\dim
Z_a>0$ or $\#Z_a\in \{0,1\}$ holds. The last statement follows by an
elementary calculation of commutators.
\end{proof}

Now we start with the

\begin{prop}
The group $G$ is not solvable.
\end{prop}
\begin{proof}
By proposition~\ref{No_C3}, assume that $G\cong G_1$, $G_2(\tau)$,
$G_3$ or $H_3$.

If $G\cong G_1$ or $H_3$, then for all $g\in G$ we have that $\dim
Z_G(g)\geq 2$. Thus, by the 2-dimensional normalizer argument, $\Gamma
< Z_G$. A direct calculation shows that $Z_G$ is contained in a
connected Abelian subgroup, contrary to proposition~\ref{No_C3}.

Now assume that $G\cong G_2(\tau)$ or $G_3$. Note that $G\cong
A\ltimes G'$ where $G'=(0,\mathbb C,\mathbb C)$ is the commutator
group of $G$ and $A$ as well as $G'$ are connected and
Abelian. Furthermore, $\dim Z_G(g')=2$ for all $g'\in G'$. Thus, by
the 2-dimensional normalizer argument $\Gamma\cap G'=\{e\}$, and the
projection $\pi_1:G\to A$ is an injective group morphism, if
restricted to $\Gamma$. This shows that $\Gamma$ is Abelian, which in
turn implies that $\Gamma<N_G(g)$.

If $\dim Z_G(\gamma)\geq 2$ for all $\gamma \in \Gamma$, then $\Gamma$
is again central and contained in a connected Abelian subgroup.  Thus
the same argument as above applies. Thus we may assume that there
exists $\gamma \in \Gamma$ with $\dim N_G(\gamma)=1$. If $N_G(\gamma)$
is connected, then $\Gamma \subset N_G(\gamma)$ and we obtain a
contradiction to lemma~\ref{Abelian_lemma}. If $N_G(\gamma)$ is not
connected, then lemma~\ref{centralizer_lemma} asserts that
$Z_G(\gamma)^0\subset G'$. But this is also not possible: use the
1-dimensional normalizer argument to see that there exists an element
$\gamma'\in \Gamma \cap G'$. However, above we have already ruled out
this possibility.
\end{proof}

This finishes the proof of the main theorem \ref{mainThm} up to the
proof \ref{CCAA:prop} and the $\mathbb C^{*}\times \mathbb
C^{*}$-action argument.

\section{Proof of the $\mathbb C^{*}\times \mathbb C^{*}$-action argument}

\label{CCAA:Section}

\begin{lem}
Suppose that $T\cong (\mathbb C^{*})^2$ acts on $X$. If $E_{i}$
is an irreducible component of $E$, then $T$ acts almost transitively
on $E_{i}$.
\end{lem}
\begin{proof}
Suppose that $T$ does not act almost transitively on $E_{i}$. Then
proposition \ref{Fix.of.C*} asserts that the generic $T$-orbit
$O\subset E_{i}$ is 1-dimensional. We know that $O\cong \mathbb C^{*}$
or that $O$ is isomorphic to a 1-dimensional compact torus. In either
case, $\dim Aut(O)=1$, so that there is a 1-dimensional kernel
$T'=ker(T\to Aut(O))$. If $x\in O$ is any point, then linearize the
$T'$-action at $x$ (note that $O$ does not meet the singular locus of
$E_{i}$). The linearization shows that either $T'$ acts trivially on
$E_{i}$, or that the $T'$-orbits are transversal to $O$. The first
case is ruled out by proposition \ref{Fix.of.C*}.  The second case is
ruled out by the assumption that $T$ does not act almost
transitively. In each case we obtain a contradiction.
\end{proof}
\begin{lem}
There exists an irreducible component $E_{i}$ of $E$ which is rational.
\end{lem}
\begin{proof}
First, we claim that the $T$-fixed point set is not empty. Identifying
$T$ with $(\mathbb C^{*})^2$, we write $T=T_1.T_2$. By
lemma \ref{Fix.of.C*} we know that $Fix(T_1)$ consists of 2 points,
or it is isomorphic to $\mathbb P_1$. In both cases we are finished if
we note that $T$ is Abelian, so that $T$ stabilizes $Fix(T_1)$. 

Second, choose a point $x$ in the $T$-fixed point set and a number $i$
such that $x\in E_{i}$. We will show that there is a $T$-stable
rational curve $C\subset E_{i}$ with $x\in C$. Then, if $\hat{E}_{i}$
is a minimal desingularization and $\hat{C}$ a component of the
preimage of $C$ which is not mapped to a point, $\hat{C}$ is still
$T$-stable, and contains a $T$-fixed point. The classification of the
smooth almost homogeneous surfaces yields that $\hat{E}_{i}$ contains
a $T$-fixed point only if it is rational. See
e.g. \cite[p. 92f]{H-Oe}.

In order to construct $C$, linearize the $T$-action at $x$. Denote the
weights of $T_1$ by $(a,b,c)$ and the weights of $T_2$ by
$(d,e,f)$. Recall from proposition \ref{Fix.of.C*} that at most one
weight in each triple is 0. If necessary, use an embedding $\mathbb
C^{*}\to T$, $\lambda \mapsto (\lambda ^{-d},\lambda ^{a})$ to find an
action of $T'\cong \mathbb C^{*}$ on $X$, fixing $x$ and having
weights $(0,g,h)$. Now let $C$ be the $T'$-fixed point curve. Since
$c_1(\mathcal{O}_X(E_i))=0$, $C$ must be contained in $E_{i}$.
\end{proof}
Now we finish the proof of the $\mathbb C^{*}\times \mathbb C^{*}$-action
argument \ref{CCAA:prop}.

\begin{proof}
Let $E_{i}\subset E$ be a rational component, and $E_{j}$ a component
different from $E_{i}$.  If $\delta :\hat{E}_{i}\to E_{i}$ is a
minimal resolution of the singularities, then, since $b_2(X)=0$,
$c_1(\delta^*(\mathcal{O}_X(E_j)))=0$, but $\delta ^{*}(E_{j})$ is an
effective divisor. This contradicts $\hat{E}_{i}$ being rational.
\end{proof}

\section{Considerations concerning the $\mathbb C^{*}$-action argument}
\label{StartC*}

The remainder of the paper is devoted to proving proposition~
\ref{S1Proposition}. As was indicated at the beginning of
section~\ref{sect:setup}, this completes the proof of our main
theorem, i.e., that $X$ is not almost homogeneous.

Accepting the situation presented to us by proposition~
\ref{S1Proposition}, we operate here under the assumptions that
$\Omega $ has the Betti type of $S^1$ and that there exists an
effective $\mathbb C^*$-action on $X$.  We begin by deriving some
topological consequences of the assumption on the Betti type of
$\Omega $.

\subsection{Topological constraints}

\begin{prop}\label{CAA:E_is_rational}
The divisor $E$ is irreducible and not normal. Its normalization
$\tilde{E}$ has only rational singularities and the minimal
desingularization $\hat{E}$ is rational. In particular, $\tilde{E}$ is
$\mathbb Q$-factorial.
\end{prop}

\begin{proof}
Since $\Omega $ has Betti-type $S^1$, it follows from
proposition~\ref{Top:EvsOmega} that $b_{4}(E)=1$, i.e. $E$ is
irreducible. Suppose that $E$ is normal and let $\pi :\hat{E}\to E$ be
the minimal desingularization. Since $b_3(E)=0$ we have
$b_3(\hat{E})=0$ (see the proof of proposition \ref{TProposition}),
hence $b_1(\hat{E})=0$, and $\hat{E}$ is K\"ahler.

Since $K_{\hat{E}}\subset \pi ^{*}(K_{E})$ and $\pi ^{*}(K_{E})$ is
linearly equivalent to 0, we have $\kappa (\hat{E})\leq 0$, and by
$b_3(\hat{E})=0$, $\hat{E}$ is either rational or birational to a
K3-surface or an Enriques surface. But because $\hat{E}$ possesses a
$\mathbb C^{*}$-action, the latter cases are excluded. So $\hat{E}$
(and hence $E$) are rational surfaces. As a consequence note that
$R^1\pi _{*}(\mathbb Q)=0$ by Leray's spectral sequence and
$$
H^1(\hat{E},\mathbb Q)=H^2(E,\mathbb Q)=0.
$$
Thus $H^2(\hat{E},\mathbb Q)=H^0(E,R^2\pi _{*}(\mathbb Q))$ and
we conclude that $H^2(\hat{E},\mathbb Q)$ is generated by the
$\pi$-exceptional curves.

Now take an ample divisor $A$ on $\hat{E}$. Then we find $m\in \mathbb
N$, $\lambda _{i}\in \mathbb Z$ and $\pi $-exceptional curves
$C_{i}\subset \hat{E}$ such that
$$
c_1(\mathcal{O}_{\hat{E}}(mA)) =
c_1\left(\mathcal{O}_{\hat{E}}\left(\sum \lambda_i Ci\right)
\right). 
$$ 
Since $\hat{E}$ is rational we have the linear equivalence $mA=\sum
\lambda _{i}C_{i}$ which is absurd.

Consequence: $E$ is not normal.

Let $\nu :\tilde{E}\to E$ be the normalization, and $\pi :\hat{E}\to
\tilde{E}$ the minimal desingularization. Using the formula $\omega
_{\tilde{E}}=\nu ^{*}(\omega _{E})-\tilde{N}$ and the fact that
$\tilde{N}$ is an effective Weil-divisor supported on the
preimage of the non-normal locus (observe that $E$ is Gorenstein!),
the ``old'' arguments still apply and give the rationality of $E$.

In order to show that $\tilde{E}$ has only rational singularities, we
check 
$$
R^1\pi_*(\mathcal{O}_{\hat{E}}) = 0.
$$ 
Since $H^1(\mathcal{O}_{\hat{E}})=0$, Leray's spectral sequence
yields an embedding 
$$
H^0(R^1\pi_*(\mathcal{O}_{\hat{E}}))\to H^2(\mathcal{O}_{\tilde{E}}).
$$ 
Now $\tilde{E}$ is Cohen-Macaulay and therefore
$H^2(\mathcal{O}_{\tilde{E}}) \cong H^0(\omega_{\tilde{E}})$. Since
$\omega _{\tilde{E}}\subset \nu ^{*}(\omega _{E})$ and
$\omega_{\tilde{E}}\not = \nu^*(\omega_E)=\mathcal{O}_E$, we have
$H^2(\mathcal{O}_{\tilde{E}}) = 0$. Therefore
$R^1\pi_*(\mathcal{O}_{\hat{E}})=0$.
\end{proof}

\begin{notation}
Let $\nu :\tilde{E}\to E$ be the normalization and $\pi :\hat{E}\to
\tilde{E}$ be the minimal desingularization. Write $\delta :=\pi \circ
\nu $. Let $N\subset E$ be the non-normal locus and $\tilde{N}:=\nu
^{-1}(N)$.
\end{notation}

\begin{prop}\label{CAA:Bettinum.of.N.tilde.N}
The following Betti-numbers are equal: $b_0(N) =
b_0(\tilde{N})$ and $b_1(N)=b_1(\tilde{N})$.
\end{prop}

\begin{proof}
We use the following Mayer-Vietoris\label{red_Mayer-Vietoris}
sequence for reduced cohomology:
$$
\ldots \to \tilde{H}^{q}(E;\mathbb Q)\to 
\tilde{H}^{q}(\tilde{E};\mathbb Q)\oplus \tilde{H}^{q}(N;\mathbb Q)\to 
\tilde{H}^{q}(\tilde{N};\mathbb Q)\to \tilde{H}^{q+1}(E;\mathbb Q)\to \ldots 
$$
(see \cite[prop. 3.A.7 on p. 98]{BK82} for information about this sequence).

So $h^0(\tilde{N};\mathbb Q)=h^0(N;\mathbb Q)$ since
$H^1(E;\mathbb Q)=0$ by proposition~\ref{Top:EvsOmega}. Furthermore
$b_1(\tilde{E})=b_1(\tilde{N})-b_1(N)$, since $H^2(E;\mathbb
Q)=0$.
\end{proof}

\begin{lem}\label{CAA:14}
The space $N$ is connected.
\end{lem}
\begin{proof}
Using the fact that $E$ is Cohen-Macaulay, by \cite[p.~166,
3.34(2)]{Mori82} there exists an exact sequence
$$
\begin{CD} 0@>>> \mathcal{O}_E @>>>
\nu_*(\mathcal{O}_{\tilde{E}}) @>>> \omega^{-1}_E\otimes\omega_N @>>>
0\end{CD}
$$ 
Since $c_1(\nu ^{*}(\omega _{E}))=0$ so that
$\nu^*(\omega_E)=\mathcal{O}_{\tilde{E}}$,
$$
\begin{CD} 0 @>>> \omega_E
@>>>
\underbrace{\omega_E\otimes\nu_*(\mathcal{O}_{\tilde{E}})}_{=\nu_*\nu^*
(\omega_E)=\nu_*(\mathcal{O}_{\tilde{E}})}
@>>> \omega_N @>>> 0\end{CD}.
$$ 
is also exact. Again using \cite[p. 166, 3.34(2)]{Mori82} we see that
$N$ is Cohen-Macaulay so that Serre duality holds. Thus,
the associated long cohomology sequence gives $$\begin{CD} \ldots @>>>
\underbrace{H^1(\tilde{E},\mathcal{O}_{\tilde{E}})}_{=0} @>>>
\underbrace{H^1(N,\omega_N)}_{\cong H^0(N,\mathcal{O}_N)} @>>>
\underbrace{H^2(E,\omega_E)}_{\cong H^0(E,\mathcal{O}_E) \cong \mathbb
C} @>>> \ldots \end{CD}$$ so that $h^0(N,\mathcal{O}_N)\leq 1$, and
the reduced subspace $N_{red}$ is connected.
\end{proof}

\subsection{Orbits of the $\mathbb C^{*}$-action}

For the sake of completeness we outline some known facts on $\mathbb
C^{*}$-actions on rational surfaces: the orbits are always
constructible. As a consequence we see that $E$ necessarily contains
attractive and repulsive fixed points, a fact that will be crucial in
the sequel. More information can be found in the works of Bia\l
ynicki-Birula and of Sommese (see e.g.~\cite{BBS85}).

\begin{lem}\label{CAA:C*_is_Closed}
If $H$ is a linear algebraic group and $\iota :\mathbb C^{*}\to H$ is
a holomorphic map which is a (set-theoretical) group morphism, then
the image $\iota (\mathbb C^{*})$ is a closed algebraic subgroup of
$H$.
\end{lem}
\begin{proof}
Let $J<H$ be the Zariski closure of $\iota (\mathbb C^{*})$, in
$H$. It is immediate that $J$ is Abelian and, since it is affine,
$J\cong \mathbb C^{n}\times (\mathbb C^{*})^{m}$.

Now use the fact that there is no non-constant holomorphic group morphism from
$\mathbb C^{*}$ to $\mathbb C$, and that any group morphism from $\mathbb C^{*}$to
$\mathbb C^{*}$ is given by $z\mapsto z^{k}$.
\end{proof}

\begin{lem}\label{CAA:1}
Let $S$ be a (possibly singular) irreducible surface with a
holomorphic action of $\mathbb C^{*}$. Suppose that the minimal
desingularization $\hat{S}$ is rational. Then:
\begin{enumerate}
\item All $\mathbb C^{*}$-orbits are constructible, and their closures 
are rational.  In particular, for all $x\in S$, the limits $\lim
_{\lambda \in \mathbb C^{*},\lambda \to 0}\lambda .x$ and $\lim
_{\lambda \in \mathbb C^{*},\lambda \to \infty }\lambda .x$ exist.
\item If $F\subset S$ is the set of the $\mathbb C^{*}$-fixed
points. Then there are two components $F_0$ and $F_{\infty }$ of $F$
and a Zariski-open set $U\subset S$ such that for all
$x\in U$:
$$
\lim _{\lambda \in \mathbb C^{*},\lambda \to 0}\lambda .x\in F_0,\quad 
\lim _{\lambda \in \mathbb C^{*},\lambda \to \infty }\lambda .x\in F_{\infty }.
$$
\end{enumerate}
\end{lem}

\begin{proof}
Suppose for the moment that $S$ was smooth. Then $Aut(S)$ is linear
algebraic and acts algebraically; this is a consequence of the fact
that, since $b_1(S)=0$, $S$ can be equivariantly embedded into some
$\mathbb P_{n}$ ---see \cite{Bla56} for a proof. By lemma
\ref{CAA:C*_is_Closed}, any closed subgroup $\mathbb C^{*}<Aut(S)$ is
linear algebraic and acts algebraically.  In particular, $\mathbb
C^{*}$-orbits are constructible. It follows from Borel's fixed
point theorem (see e.g. \cite[p.~32]{H-Oe}) that all $\mathbb
C^{*}$-stable curves in $S$ contain fixed points. This already proves
(1).

In order to prove assertion (2), embed $\mathbb C^{*}\to \mathbb P_1$
in the usual way. Then there exists a rational morphism 
$$
\phi :\mathbb P_1\times S\dasharrow S.
$$
Since the set of fundamental points of $\phi ^{-1}$ is of codimension
$\geq 2$, there exists an open set $U\subset S$ such that $\phi
|_{\mathbb P_1\times U}$ is regular. Consequence: for all $x\in U$
the limits $\lim _{\lambda \in \mathbb C^{*},\lambda \to 0}\lambda .x$
and $\lim _{\lambda \in \mathbb C^{*},\lambda \to \infty }\lambda .x$
exist.  Set $F_0=\overline{\phi (0\times U)}$ and $F_{\infty
}:=\overline{\phi (\infty \times U)}$.

If $S$ is singular, then let $\delta :\hat{S}\to S$ be an equivariant
resolution and find $\hat{F}_0$, $\hat{F}_{\infty }$ and
$\hat{U}\subset \hat{S}$ as above. It is sufficient to set
$F_0:=\delta (\hat{F}_0)$, $F_{\infty }:=\delta (\hat{F}_{\infty
})$ and $U:=\delta (U)\cap S_{reg}$, where $S_{reg}$ is the (open) set
of regular points in $S$.
\end{proof}

\section{The case of a discrete fixed point set}\label{sect:discrete_fixed_points}

Our goal here is to prove the $\mathbb C^*$-action argument under the
assumption that the set $F$ of $\mathbb C^*$-fixed points is discrete,
i.e.~$F=\{F_0, F_\infty\}$. It is shown that the discreteness
assumption, which is made throughout this section, is contrary to
$b_1(N)=b_1(\tilde N)$ (see proposition~\ref{CAA:no_discrete}).

\begin{notation}
In the sequel $U_0$ and $U_\infty$ are disjoint linearizing
neighborhoods of the points $F_0$ and $F_\infty $ in $X$.
\end{notation}

\subsection{Symmetry lemmas}

Here we prove several lemmas which show that the situations at $F_0$
and at $F_{\infty }$ are very similar. First, we investigate:

\subsubsection{Weights of the $\mathbb C^{*}$-actions on $T_{F_0}X$
and $T_{F_{\infty }}X$ }

The main result of this section is

\begin{prop}
The locally closed spaces $E\cap U_0$ and $E\cap U_{\infty }$ are
reducible. One of the components of each space is smooth, and the
$\mathbb C^{*}$-action has a totally attractive resp. repulsive fixed
point there.
\end{prop}
Since the proof is somewhat lengthy, and we will use a number of
partial results later, we subdivide the proof into a sequence of
lemmas and corollaries.

\begin{lem}\label{CAA:3}
After swapping $F_0$ and $F_\infty$ and, if necessary, replacing the
$\mathbb C^{*}$-action by $(\mathbb C^{*})^{-1}$, the weights of the
$\mathbb C^{*}$-action on the tangent space $T_{F_0}X$ have the
following signs: $(++-)$.
\end{lem}
\begin{proof}
By symmetry, we only have to exclude the following distribution of
signs of the weights of the $\mathbb C^{*}$-action on $T_{F_0}X$ and
$T_{F_{\infty }}X$:
\begin{description}
\item [$\boldsymbol{(+++),(+++)}$]This contradicts lemma \ref{CAA:1}. 
\item [$\boldsymbol{(+++),(---)}$]Let $x\in E$ be a generic point. Choose 
a sequence $(x_{n})_{n\in \mathbb N}\subset \Omega $ such that $\lim
_{n\to \infty }x_{n}=x$.  Then there exist numbers $\lambda $,
$\lambda '\in \mathbb C^{*}$ such that $\lambda .x\in U_0$ and
$\lambda '.x\in U_{\infty }$ and there exists an $n\in \mathbb N$ such
that $\lambda .x_{n}\in U_0$, and $\lambda 'x_{n}\in U_{\infty
}$. Linearization shows that $\lim _{\lambda \to 0}\lambda
x_{n}=F_0$ and $\lim _{\lambda \to \infty }\lambda x_{n}=F_{\infty
}$ so that $C:=\overline{\mathbb C^{*}.x_{n}}$ is a closed curve in X
and $C\cap E=F_0\cup F_{\infty }$. But then $E.C\geq 2$, a
contradiction to $H^2(X,\mathbb Z)=0$.
\end{description}
\end{proof}
Assume from now on that the $\mathbb C^{*}$-action on $T_{F_0}X$ has
weights of type $(++-)$. Let $x$, $y$ and $z$ be coordinates for the
associated weight spaces. By the linearization theorem~
\ref{linearization}, we can view $x$, $y$ and $z$ as giving local
coordinates on $U_0$.  After performing a $\mathbb C^{*}$-equivariant
change of coordinates on $T_{F_0}X$, we may assume that the unit ball
in $T_{F_0}X$ is contained in the image of $U_0$.

\begin{cor}\label{CAA:4} 
After swapping $F_0$ and $F_{\infty }$ and, if necessary, replacing
the $\mathbb C^{*}$-action by $(\mathbb C^{*})^{-1}$, the weights of
the $\mathbb C^{*}$-action on the tangent space $T_{F_0}X$ have signs
$(++-)$ and the locally closed subspace $E\cap U_0$ is reducible. One
of the components of $E\cap U_0$ is smooth, and the $\mathbb
C^{*}$-action has a totally attracting fixed point there.
\end{cor}
\begin{proof}
Since $E\cap U_0$ is closed in $U_0$ and contains infinitely many
$\mathbb C^{*}$-stable curves containing $F_0$, we have that
$\{z=0\}\cap U_0\subset E$.

Since $\{z=0\}$ is smooth, in order to see that $E\cap U_0$ is
reducible, it suffices to show $E\cap U_0$ is not normal. Since $E$ is
a hypersurface in $X$, it is Cohen-Macaulay, and it follows
from Serre's criterion that the non-normal locus $N\subset E$
must be of codimension 1.  By lemma \ref{CAA:1}, $N$ contains $\mathbb
C^{*}$-fixed points.  If $N$ contains $F_0$, then we can stop here.

Otherwise, if $N$ contains $F_{\infty }$, but not $F_0$, then
the signs of the weights at $F_{\infty }$ are necessarily $(--+)$;
we show this by ruling out all other possibilities:
\begin{description}
\item [$\boldsymbol{(+++)}$]does not occur, or else obtain a contradiction 
to lemma \ref{CAA:1}. 
\item [$\boldsymbol{(++-)}$]similarly
\item [$\boldsymbol{(---)}$]then $\lim _{\lambda \to 0}\lambda x$ would 
not exist for generic $x\in N$. Since every 1-dimensional component of
$N$ contains a $\mathbb C^{*}$-fixed point, it is rational. Thus, the
limit exists on the normalization of $N$ which in turn implies that
the limit exists on $N$.
\end{description}
In this case swap $F_0$ and $F_{\infty }$ and start anew.
\end{proof}
For the remainder of this section, fix $F_0$ and $F_\infty$ so that we
are in the situation of the above corollary.

\begin{notation}
Let $E_{0,i}$ denote the irreducible components of $E\cap U_0$ and let
$E_{0,0}$ be the smooth component in the preceding corollary.
\end{notation}
\begin{lem}\label{CAA:5}
Choose numbers $a$, $b$ and $c\in \mathbb N^{+}$and let the group
$\mathbb C^{*}$ act on $\mathbb C^3$ by $\lambda :(x,y,z)\to (\lambda
^{a}x,\lambda ^{b}y,\lambda ^{-c}z)$.  Let $S\subset \mathbb C^3$ be
an irreducible $\mathbb C^{*}$-stable divisor with $S\not =\{z=0\}$,
but $S\cap \{z=0\}\not =\emptyset $.  Then $\{x=y=0\}\subset S$.
\end{lem}
\begin{proof}
Choose a point $s\in \{z=0\}\cap S$, s. Let
$s_{n}=(x_{n},y_{n},z_{n})$ be a sequence with $z_{n}\not =0$ and
$\lim s_{n}=s$. Choose $\lambda _{n}\in (z_{n})^{\frac1{c}}$.  Note
that $\lim \lambda _{n}=0$. Then $\lambda _{n}.s_{n}\in E$, and $\lim
\lambda _{n}s_{n}=\lim (\lambda ^{a}_{n}x_{n},\lambda
^{b}_{n}y_{n},1)=(0,0,1)$.
\end{proof}

\begin{cor}\label{CAA:6}
Every component $E_{0,i}$ $(i\not =0)$ contains $\{x=y=0\}$.
\end{cor}

\begin{lem}
\label{CAA:7}If the signs of the weights of the $\mathbb C^{*}$-action
on $T_{F_{\infty }}X$ are all negative, then there exists a curve $C^{-}\subset E$
satisfying
\begin{enumerate}
\item $C^{-}\cap E_{0,0}$ is a curve (i.e. $\dim C^{-}\cap E_{0,0}=1$) and
\item $F_{\infty }\not \subset C^{-}$
\end{enumerate}
\end{lem}
\begin{proof}
By corollary \ref{CAA:4}, $E\cap U_0$ is reducible, and by corollary
\ref{CAA:6}, $\{x=y=0\}\subset E\cap U_0$. The $z$-axis is the
weight-space to the negative weight, so that $\lim _{\lambda \to
\infty }(0,0,1)=F_0$.  But the limit $\lim _{\lambda \to 0}(0,0,1)$
exists; this is because $\hat{E}$ is rational, and the limit exists
there. Due to the negative weights it is impossible that $\lim
_{\lambda \to 0}(0,0,1)=F_{\infty }$. Thus $\lim _{\lambda \to
0}(0,0,1)=F_0$ and for all $\lambda $ sufficiently small $\lambda
.(0,0,1)\in E_{0,0}$.
\end{proof}

\begin{lem}\label{CAA:8}
The weights of the $\mathbb C^{*}$-action on $T_{F_{\infty }}X$ have
signs $(--+)$.
\end{lem}
\begin{proof}
It is clear that at least two of the signs must be negative ---this is
because for generic $x\in E$, $\lim _{\lambda \to \infty }\lambda
.x=F_0$.  Now suppose the weights were $(a,b,c)$ which were all
negative. The weights of the $\mathbb C^{*}$-action on
$T_{F_0}E_{0,0}$ shall be denoted by $d$ and $e$. We consider the
weighted projective spaces $\mathbb P_{(-a,-b,-c)}(T_{F_{\infty }}X)$
and $\mathbb P_{(d,e)}(T_{F_0}E_{0,0})$. These are parameter spaces
for $\mathbb C^{*}$-stable curves in $X$ and $E$ passing through $F_0$
or $F_{\infty }$, respectively.

The analytic subspace $E\cap U_{\infty }$ gives a \emph{closed}
subspace in the weighted projective space $\mathbb E\subset \mathbb
P_{(-a,-b,-c)}(T_{F_{\infty }}X)$ parameterizing curves in $E\cap
U_{\infty }$. We will now construct a map from $\mathbb E$ to $\mathbb
P_{(d,e)}(T_{F_0}E_{0,0})$.

First fix some notation: let $\Lambda _0:E_{0,0}\to T_{F_0}E_{0,0}$
and $\Lambda _{\infty }:U_{\infty }\to T_{F_{\infty }}X$ be the
linearizing maps, and let $\pi _0:T_{F_0}E_{0,0}\setminus \{0\}\to
\mathbb P_{(d,e)}(T_{F_0}E_{0,0})$ be the canonical projection.

Given an arbitrary point $x\in \mathbb E$, there exists a neighborhood
$U(x)\subset \mathbb E$ and a section $\sigma :U(x)\to \Lambda
_{\infty }(U_{\infty })\subset T_{F_{\infty }}X$.  After shrinking
$U(x)$, if necessary, there is a $\lambda \in \mathbb C^{*}$ such that
$(\lambda \circ \sigma )(U(x))\subset E_{0,0}\setminus \{F_0\}$:
simply choose a $\lambda $ such that $(\lambda \circ \sigma )(v)\in
E_{\infty }$ and set $U'(x):=\lambda ^{-1}((\lambda \circ \sigma
)(U(x))\cap E_{0,0})$; this is the shrinkage that might be
unavoidable. This way we obtain a map $\iota _{x}:=(\pi _0\circ
\lambda \circ \sigma ):U(x)\to \mathbb P_{(d,e)}(T_{F_0}E)$.

In order to obtain a global map $\iota :\mathbb E\to \mathbb
P_{(d,e)}(T_{F_0}E)$, it suffices to show that $\iota _{x}$ does not
depend on the choice of $\sigma $ and $\lambda $. Indeed, choosing
different $\sigma '$ and $\lambda '$ then for all $y\in U(x)$ there is
a unique number $\lambda _{y}\in \mathbb C^{*}$ such that $(\lambda
\circ \sigma )(y)=\lambda _{y}(\lambda '\circ \sigma ')(y)$.  This
already shows that $(\pi _0\circ \lambda \circ \sigma )(y)=(\pi
_0\circ \lambda '\circ \sigma ')(y)$, and the existence of the global
map $\iota $ is shown. It is obvious that $\iota $ is injective.

This is how we make use of $\iota $ : by lemma \ref{CAA:7}, $\pi
(C^{-}\cap E_{0,0})\setminus \{F_0\}$ is not contained in the image of
$\iota $, so that the image must be contained in $\mathbb
P_{(d,e)}(T_{F_0}E)\setminus (\pi (C^{-}\cap E_{0,0}\setminus
\{F_0\}))\cong \mathbb C$.  By the maximum principle, the image must
be a point. A contradiction to $\iota $ being injective.
\end{proof}

\begin{cor}
\label{CAA:9}$E\cap U_{\infty }$ is reducible. There exists a smooth component
$E_{\infty ,0}$ with totally repulsive fixed point.
\end{cor}
\begin{proof}
The existence of $E_{\infty ,0}$ follows exactly as in corollary \ref{CAA:4}.
Similarly, it follows from the argumentation of corollary \ref{CAA:4} that
$E\cap U_{\infty }$ is reducible, if we show that the non-normal locus
$N\subset E$ intersects $E_{\infty ,0}$. 

Suppose this was not the case. Then, if $x\in E_{\infty ,0}$, $\lim _{\lambda \to 0}\lambda x=F_0$,
and the argumentation used in the proof of lemma \ref{CAA:8} yields a contradiction.
\end{proof}
\begin{notation}
\label{CAA:10}In analogy to the notation introduced above, let $E_{\infty ,i}$
denote the irreducible components of $E\cap U_{\infty }$ and let $E_{\infty ,0}$
be the smooth component whose existence is asserted by the preceding corollary.
\end{notation}

\subsubsection{Loops}

Now turn to the non-normal locus of $N\subset E$ . The following is an
important notion:

\begin{notation}
\label{CAA:loop.notation}Call a $\mathbb C^{*}$-stable curve $C\subset E$
a ``loop'', if for a generic point $x\in C$:
$$
\lim _{\lambda \to 0}\lambda .x=\lim _{\lambda \to \infty }\lambda .x$$
 
\end{notation}
\begin{lem}
\label{CAA:11}\label{CAA:12}There is at most one loop containing $F_0$,
and at most one containing $F_{\infty }$. The number of loops containing
$F_0$ and the number of loop containing $F_{\infty }$ are equal.
\end{lem}
\begin{proof}
There is only one curve in $E\cap U_0$ containing a point $x$ such
that $\lim _{\lambda \to \infty }\lambda .x=F_0$ (namely the $z$-axis).
The situation at $F_{\infty }$ is similar. 

Argue as in the proof of lemma \ref{CAA:8}, using a map 
$$
\mathbb P_{(a,b)}(T_{F_0}E_{0,0})\to \mathbb P_{(c,d)}(T_{F_{\infty }}E_{\infty ,0})$$
to exclude the possibility that there is no loop at $F_0$ and one at
$F_{\infty }$ or vice versa.
\end{proof}

\begin{lem}\label{CAA:16}
The number of irreducible components of $N\cap E_{0,0}$ equals the
number of irreducible components of $N\cap E_{\infty ,0}$.
\end{lem}
\begin{proof}
Decompose $N=N_{L}\cup \bigcup _{i}N_{i}$, where the $N_{L}$ are loops
and the $N_{i}$ are other components. By lemma \ref{CAA:11}, the claim
is true if one considers loops only. If $N_{i}$ is one of the other
components and $x\in N_{i}$ a generic point, then $\lim _{\lambda \to
0}\lambda .x=F_0$ and $\lim _{\lambda \to \infty }\lambda .x=F_{\infty
}$ so that $N_{i}\cap E_{0,0}$ and $N_{i}\cap E_{\infty ,0}$ are both
irreducible components of $N\cap E_{0,0}$ and $N\cap E_{\infty ,0}$,
respectively.
\end{proof}

\subsection{Preparations: $\mathbb C^*$-action on normal surfaces}

Let $S$ be a smooth connected algebraic surface equipped with an
algebraic $\mathbb C^*$-action with an attractive (resp. repulsive)
fixed point $F_\infty$ (resp. $F_0$). Let $U\subset S$ be as in
lemma~\ref{CAA:1}.

\begin{notation}
A curve $C=C_+\cup C_1\cup \ldots\cup C_k\cup C_-$, $k\geq 0$, as in
the following figure
$$
\xymatrix {
\ _{F_0} \bullet \ar@{-}[rr]|{>}_{C_+}& &
 \bullet_{p_1} \ar@{-}[rr]|{>}_{C_1}& &
 \bullet_{p_2} \ar@{-}[rr]|{>}_{C_2}& &
 \cdots  \ar@{-}[rr]|{>}_{C_-}& & 
 \bullet_{F_\infty}
}
$$
which is invariant under the $\mathbb
C^*$-action will be called an ``external chain''.

A $\mathbb C^*$-fixed point different from $F_0$ and $F_\infty$ is an
``external fixed point''.
\end{notation}

\begin{lem}\label{CAA:complement}
The complement of $U$ in $S$ is a union of external chains without
common components.
\end{lem}
\begin{proof}
At every external fixed point the weights of the linearization must be
of type $(+-)$.
\end{proof}

Now allow $S$ to have normal singularities at external fixed
points. Applying the above argument to the desingularization $\hat S$
and blowing back down, we have the same result.

\begin{lem}
If $S$ is a normal compact rational surface equipped with an $\mathbb
C^*$-action having a smooth point $F_0$ as a source and a smooth point
$F_\infty$ as a sink, then the complement $S\setminus U$ is a union of
external chains.
\end{lem}

\begin{cor}
In the setting of the preceding lemma there are no loops in $S$.
\end{cor}

\subsection{Computation of $b_1$ of curves}

\begin{lem}\label{norb1}
Let $C$ be the union of $c$ rational curves having $s$ singular
points. Let $\nu:\tilde C\to C$ be the normalization of $C$ and
$\tilde s$ be the number of points in $\nu^{-1}(C_{Sing})$. For
convenience, let $\delta=\tilde s-s$. Then
$$
b_1(C)=1-c+\delta
$$
\end{lem}
\begin{proof}
Apply the Mayer-Vietoris sequence (see
p.~\pageref{red_Mayer-Vietoris}), where $E:=C$, $N:=C_{Sing}$, etc.:
$$
0  = \underbrace{\tilde h^0(E)}_{0} -  (\underbrace{\tilde
h^0(\tilde E)}_{c-1} + \underbrace{\tilde h^0(N)}_{s-1}) + \underbrace{\tilde
h^0(\tilde N)}_{\tilde s-1}-\underbrace{h^1(E)}_{b_1}+(\underbrace{h^1(\tilde E)+h^1(N)}_{0}).
$$
\end{proof}

\subsection{The proof in the case where $F$ is discrete}

\begin{notation}
An ``outer orbit'' is an orbit which flows in the opposite direction
from the generic orbit, i.e.~ with source $F_\infty$ and sink $F_0$.

An irreducible $\mathbb C^*$-invariant curve $C$ containing $F_0$ as a
source and $F_\infty$ as a sink is called a ``crossing curve'' if $E$
is not locally irreducible at the generic point of $C$. 
\end{notation}

\begin{prop}\label{CAA:no_discrete}
The case where $F$ is discrete does not occur.
\end{prop}
\begin{proof}
Assume to the contrary and use the notation introduced above. Let $n$
be the number of components in $N\cap E_{0,0}$. Then, by
lemma~\ref{norb1},
$$
b_1(N)=\left\{
\begin{matrix}
n-1 & \text{ if $N$ consists of crossing curves only} \hfill \\ 
n & \text{ if $N$ contains an outer orbit} \hfill \\
n+1 & \text{ if $N$ contains loops } \hfill
\end{matrix}
 \right.
$$ 
on the other hand, since the number of external chains in $\tilde E$
plus the number of curves in $\tilde U\cap \tilde N$ is also $n$ (note
that there are $n$ components of $\tilde N$ containing $\tilde F_0$),
we have
$$
b_1(\tilde N)=n-1.
$$

It remains to show that the case where $N$ only consists
of crossing curves does not occur.  In that case, since all
components $E_{0,i}$, $i\ne 0$, contain $\{ x=y=0\} $ and 
$N$ does not contain an outer orbit, there is only one such
component $E_{0,1}$.  Furthermore, $E_{0,0}\cap E_{0,1}$
consists of closures of $\mathbb C^*$-orbits which flow from
$F_0$ to $F_\infty $.

Observe that in this situation there can not be external chains in
$\tilde E$: If $p\in \tilde N$ is an external fixed point, then it is
the intersection point of a curve flowing {\em into} $p$ with a curve
flowing {\em out of} $p$. But if $\nu(p)=F_0$, then there all curves
in $N=\nu (\tilde N)$ flow {\em out of} $\nu(p)$, and we would have
contradiction.

The analogue holds if $\nu(p)=F_\infty$.
\end{proof}

\section{Proof of the $\mathbb C^{*}$-action argument if $F$
is a curve}\label{sect:fixed_curve}

\label{EndC*}If $F$ is a curve, then $F\cong \mathbb P_1$ by proposition
\ref{Fix.of.C*}. It is very easy to calculate $b_1(N)$. However, we
first have to show that $F\subset N$.

\begin{lem}
\label{CAA:loc@x}Suppose $\dim F=1$. If $x\in F$ is an arbitrary
point, and $U$ a linearizing neighborhood of $\mathbb C^{*}$ about
$x$, then the weights of the $\mathbb C^{*}$-action on $T_{x}X$
have signs $(0+-)$, and if $x$, $y$ and $z$ are coordinates
associated to the weight spaces, then $\{yz=0\}\subset E\cap U$.
\end{lem}
Recall from the theorem on linearization \ref{linearization} that $x$,
$y$ and $z$ can be viewed to give coordinates on $U$.

\begin{proof}
We have to exclude the possibility that the signs of the non-zero weights are
equal. Suppose they were both positive. Then for a point $y\in (E\setminus F)\cap U$,
the limit $\lim _{\lambda \to \infty }\lambda y$ would not exist in $F.$
This contradicts the fact that $F$ is the full $\mathbb C^{*}$-fixed
point set.

Now it is a direct consequence of lemma \ref{CAA:1} that $\{zy=0\}\subset E\cap U$.
Note that $F_0=F_{\infty }=F$.
\end{proof}
\begin{lem}
\label{CAA:b1N}If $F\cong \mathbb P_1$, then $b_1(N)=($\# of
irreducible components of $N)-1$.
\end{lem}
\begin{proof}
Decompose $N=F\cup \bigcup _{i}N_{i}$. We show by induction that
$b_1(F\cup \bigcup _{i=1\ldots k}N_{i})=k$.
\begin{description}
\item [Start, $\boldsymbol{k=0}$] It is clear that 
$b_1(F)=b_1(\mathbb P_1)=0$.
\item [Step]Since for all $x\in N_{i}$, the limits 
$\lim _{\lambda \to 0}\lambda x$ and $\lim _{\lambda \to \infty
}\lambda x$ exist and are $\mathbb C^{*}$-fixed, i.e. contained in
$F$, there are only two possibilities:

\begin{enumerate}
\item $\lim _{\lambda \to 0}\lambda x=\lim _{\lambda \to \infty }\lambda x$ for
all $x\in N_{k}$, i.e. $N_{k}$ is a loop (see the notation \ref{CAA:loop.notation}).
Then $b_1(N_{k})=1$ and $N_{k}\cap F$ is a single point.
\item $\lim _{\lambda \to 0}\lambda x\not =\lim _{\lambda \to \infty }\lambda x$.
Then the normalization $\tilde{N}_{k}\to N_{k}$ is injective and thus $b_1(N_{i})=0$.
Furthermore $N_{k}\cap F$ are two points.
\end{enumerate}
In any case the Mayer-Vietoris sequence associated to the
decomposition $F\cup \bigcup _{i=1\ldots k}N_{i}=(F\cup \bigcup
_{i=1\ldots k-1}N_{i})\cup N_{k}$ shows directly that 
$$ 
b_1(F\cup \bigcup _{i=1\ldots k}N_{i})=b_1(F\cup \bigcup _{i=1\ldots
k-1}N_{i})+1=k.  
$$
\end{description}
\end{proof}
Again we use a description of the rational surface $\tilde{E}$ to calculate
$b_1(\tilde{N})$ and finally to derive a contradiction.

\begin{lem}
\label{CAA:Ex.map}There exists a surjective morphism with connected fibers
$\phi :\tilde{E}\to \mathbb P_1$ such that the induced $\mathbb C^{*}$-action
on $\mathbb P_1$ is trivial.
\end{lem}
\begin{proof}
Let $C$ and $C'$ be two generic irreducible $\mathbb C^{*}$-stable
curves in $\tilde{E}$. Since the $\mathbb C^{*}$-fixed point set in
$\tilde{E}$ does not contain a totally attractive fixed point (this is
because there is non in $E$), $C$ and $C'$ are disjoint and do
not intersect the singular locus of $\tilde{E}$.

Since $\tilde{E}$ is rational, $C$ and $C'$ are linearly equivalent
as divisors and yield the desired map.
\end{proof}
\begin{lem}
The situation that $F\cong \mathbb P_1$ does not occur. 
\end{lem}
This finishes the proof of the $\mathbb C^{*}$-action argument \ref{S1Proposition}.

\begin{proof}
We show that $b_1(\tilde{N})=($ \# of irreducible components of
$N)-2$, contradicting lemma \ref{CAA:b1N}. In accordance with the
notation and the results of lemma \ref{CAA:1}, let $F_0$ and
$F_{\infty }\subset \tilde{E}$ denote the $\mathbb C^{*}$-fixed point
curves. These are sections for $\phi :\tilde{E}\to \mathbb P_1$, and
contained in the preimage $\nu ^{-1}(F)$. In particular, $F_0\cup
F_{\infty }\subset \tilde{N}$.

Again decompose $N=F\cup \bigcup _{i}N_{i}$. Set $\tilde{N}_{i}:=\nu
^{-1}(N_{i})$.

First, we show that $\tilde{N}_{i}\cap F_0$ and $\tilde{N}_{i}\cap
F_{\infty }$ are both single points. Realize from the description of
lemma \ref{CAA:loc@x} that $\nu |_{F_0}:F_0\to N$ is injective, and
note that if $y\in N_{i}$ is a generic point, then
$$
\tilde{N}_{i}\cap F_0\subset \nu ^{-1}(\lim _{\lambda \to 0}\lambda y)\cap F_0.
$$
The expression on the right hand side denotes a single point. A
similar argumentation holds for $F_{\infty }$.

Second, we claim that $b_1(\tilde{N}_{i})=0$. Recall that
$\tilde{N}_{i}$ are $\mathbb C^{*}$-stable, and do not contain $F_0$
or $F_{\infty }$ as an irreducible component. Thus, the irreducible
components of $\tilde{N}_{i}$ are irreducible components of $\phi
$-fibers.

This is sufficient information to apply the Mayer-Vietoris sequence.
For brevity, let $k:=(\#$ of irreducible components of $N)$. The
beginning of the cohomological Mayer-Vietoris sequence
associated to the decomposition $\tilde{N}=(F_0\cup F_{\infty })\cup
(\bigcup _{i}\tilde{N}_{i})$ yields:
$$
\underbrace{h^0(\tilde{N})}_{=1}-(\underbrace{h^0(F_0\cup F_{\infty })}_{=2}
+\underbrace{h^0(\bigcup _{i}\tilde{N}_{i})}_{\geq h^0(\bigcup _{i}N_{i})=k-1})
+\underbrace{h^0((F_0\cup F_{\infty })\cap 
(\bigcup _{i}\tilde{N}_{i}))}_{=2(k-1)}-h^1(\tilde{N})=0.
$$ 
In other words, $h^1(\tilde{N})\leq k-2$.
\end{proof}

\section{Some further remarks}

\begin{prop}\label{Intro:Prop1}
There is no non-zero three-form on $X$, i.e. $H^0(K_{X})=0$.
\end{prop}
\begin{proof}
Let $\sigma \in H^0(\Omega _{X}^3)$. Then $d\sigma =0$. Since
$H^3(X;\mathbb Q)=0$, de Rham's theorem exhibits a ${\mathcal
C}^{\infty}-$form $\eta $ of degree 2, such that $\sigma =d\eta
$. Thus by Stokes
$$
\int _{X}\sigma \wedge \overline{\sigma }=\int _{X}d\eta \wedge
\overline{\sigma }=\int _{X}d(\eta \wedge \overline{\sigma })-\int
_{X}\eta \wedge d\overline{\sigma }=0.
$$
Hence $\sigma =0$.
\end{proof}

\begin{cor} \ 
\begin{enumerate}
\item $H^3(X,{\mathcal O}_X) = 0;$ 
\item $H^0(K_X \otimes L) = 0$ for generic $L \in {\rm Pic}(X);$ 
\item $h^1(X,{\mathcal O}_X) = h^2(X,{\mathcal O}_X) + 1;$ 
\item ${\rm Pic}(X) = {\rm Pic}^0(X) \simeq H^1(X,{\mathcal O}_X)$ 
is a positive dimensional complex vector space. 
\end{enumerate}
\end{cor}

Observe that it is not at all clear that $\kappa (X)=-\infty $. The
proof of this fact seems to be as complicated as to show that there
are no divisors on $X$ at all (notice that $K_{X}$ cannot be torsion
by the above results!). See proposition \ref{Intro:Prop4} below.

We now prove the analogous statement to our main theorem for the
cotangent bundle.

\begin{prop}
We have $h^0(\Omega _{X}^1)\leq 2$.
\end{prop}
\begin{proof}
Since $X$ has no meromorphic functions, we have $h^0(\Omega^1_X) \leq
3,$ moreover if $h^0(\Omega^1_X) = 3,$ then any basis of 1-forms
$\omega_i$ is linearly independent at the general point: $$ \omega_1
\wedge \omega_2 \wedge \omega_3 \ne 0$$ as a section of $K_X =
\Omega^3_X.$ This contradicts proposition \ref{Intro:Prop1}.
\end{proof}
Again we can conclude that $h^0(\Omega _{X}^1\otimes L)\leq 2$ for
the general line bundle $L$. Finally we can prove statement (B) of the
introduction under some further assumptions:

\begin{prop}\label{Intro:Prop4} \ 
\begin{enumerate}
\item Assume that $X$ has no compact curves. Then 
$H^0(\Omega _{X}^1\otimes L)=0$ for all line bundles $L$, and 
$h^0(TX\otimes L)\leq 1$. 
\item Assume that $X$ has no divisors. Then 
$h^0(\Omega _{X}^1\otimes L)\leq 2$ for all line bundles $L$.
\end{enumerate}
\end{prop}

\begin{proof}
(1) Let $\omega \in H^0(\Omega _{X}^1\otimes L)$. Suppose $\omega \not
=0$.  Then there is a divisor $D$ (possibly empty) such that the
section $\omega \in H^0(\Omega^1_X \otimes L \otimes {\mathcal
O}_X(-D)) $ has no zeroes in codimension 1, hence has only finitely
many zeroes by our assumption.  Thus $$ c_3(\Omega^1_X \otimes L
\otimes {\mathcal O}_X(D)) = c_3(\Omega^1_X) \geq 0.$$ But $c_3(X) = -
c_3(\Omega^1_X) = 2,$ contradiction. The inequality $h^0(TX\otimes
L)\leq 1$ is an immediate consequence from $H^0(\Omega _{X}^1\otimes
L)=0$.

(2) Assume $h^0(\Omega^1_X \otimes L) = 3$ and take a basis $(\omega
_{i})$.  Then $\omega_1 \wedge \omega_2 \wedge \omega_3 $ is a
non-zero section of $K_{X}\otimes (3L)$.  Since there are no divisors
on $X$, the $\omega _{i}$ are linearly independent everywhere, hence
$\Omega^1_X \otimes L \simeq {\mathcal O}_X^{\oplus 3}.$ This
contradicts again $c_3(X)=2$.
\end{proof}


\begin{thebibliography}{CDP98}
\bibitem[BK82]{BK82}
G.~Barthel and L.~Kaup.
\newblock {\em Sur la Topologie des Surfaces complexes compactes}, chapter in
  Topologie des surfaces complexes compactes singuli\`eres, pages 61--297.
\newblock Number~80 in Semin. Math. Super. Les Presses de l'universit\'e de
  Montr\'eal, 1982.

\bibitem[BBS85]{BBS85}
A.~Bia\l ynicki-Birula and A.~Sommese
\newblock Quotients by ${\mathbb C}\sp*\times {\mathbb C}\sp*$ actions.
\newblock {\em Trans. Am. Math. Soc.}, 289:519--543, 1985.

\bibitem[Bla56]{Bla56}
A.~Blanchard.
\newblock Sur les vari{\'e}t{\'e}s analytiques complexes.
\newblock {\em Ann. Sci. Ec. Norm. Sup.}, 73:157--202, 1956.

\bibitem[Bre72]{Br72}
G.~Bredon.
\newblock {\em Introduction to Compact Transformation Groups}.
\newblock Academic Press, 1972.

\bibitem[CDP98]{CDP}
F.~Campana, J.-P. Demailly, and T.~Peternell.
\newblock The algebraic dimension of compact complex threefolds with vanishing
  second betti number.
\newblock {\em Compos. Math.}, 112:77--91, 1998.

\bibitem[ES86]{ES86}
J.~Erdman-Snow.
\newblock On the classification of solv-manifolds in dimension 2 and 3.
\newblock In {\em Journ{\'e}es Complexes Nancy}, number~10 in Revue de
  l'Institut Elie Cartan, 1986.

\bibitem[FF79]{FF79}
G.~Fischer and O.~Forster.
\newblock {Ein Endlichkeitssatz f\"ur Hyperfl\"achen auf kompakten komplexen
  R\"aumen}.
\newblock {\em J. reine angew. Mathematik}, 306:88--93, 1979.

\bibitem[HO80]{H-Oe}
A.~Huckleberry and E.~Oeljeklaus.
\newblock {\em Classification Theorems for Almost Homogeneous Spaces}.
\newblock Institut Elie Cartan, 1980.

\bibitem[Huc90]{Huc90}
A.~Huckleberry.
\newblock Actions of groups of holomorphic transformations.
\newblock In Barth, Narasimhan, and Gamkrelidze, editors, {\em Several complex
  variables VI. Complex manifolds}, number~69 in Encycl. Math. Sci., pages
  143--196. Springer-Verlag, 1990.

\bibitem[Jac62]{Jac62}
N.~Jacobson.
\newblock {\em Lie Algebras}.
\newblock Interscience Publishers, 1962.

\bibitem[KP85]{KP85}
H.~Kraft and V. L. Popov.
\newblock Semisimple group actions on the three dimensional affine space are linear. (English)
\newblock {\em Comment. Math. Helv.}, 60:466--479, 1985

\bibitem[Kra75]{Kras75}
V.A. Krasnov.
\newblock Compact complex manifolds without meromorphic functions.
\newblock {\em Math. Notes}, 17:69--71, 1975.
\newblock translation from Mat. Zametki, 17, 119-122 (1975). [ISSN 0001-4346].

\bibitem[Mor82]{Mori82}
S.~Mori.
\newblock {Threefolds whose canonical bundles are not numerically effective}.
\newblock {\em Ann. of Math.}, 116, 1982.

\bibitem[MU83]{Mu83}
S.~Mukai and H.~Umemura.
\newblock {Minimal Rational Threefolds}.
\newblock In M.~Raynaud and T.~Shioda, editors, {\em Algebraic geometry, Proc.
  Jap.-Fr. Conf., Tokyo and Kyoto 1982}, volume 1016 of {\em Lect. Notes Math},
  pages 490--518. Springer, 1983.

\bibitem[Pot69]{Pot69}
J.~Potters.
\newblock On almost homogeneous compact complex surfaces.
\newblock {\em Inv. Math. 8}, 1969.

\end{thebibliography}
\end{document}